\documentclass{amsart}

\usepackage{amsmath}
\usepackage{amsfonts}
\usepackage{pstricks}
\usepackage{mathrsfs}
\usepackage{latexsym}
\usepackage{times}

\theoremstyle{definition}

\newtheorem{theorem}{Theorem}

\DeclareMathAlphabet{\mathgothic}{U}{euf}{m}{n}
\DeclareMathAlphabet{\ams}{U}{msb}{m}{n}

\newpsobject{showgrid}{psgrid}{subgriddiv=1,griddots=5,gridlabels=6pt}

\def\psl{\mbox{PSL}}

\def\TT{\mathscr T}

\def\SS{\mathgothic{S}}

\def\l{\mathgothic{l}}
\def\s{\mathgothic{s}}

\def\o{\mathgothic{o}}

\def\ve{\varepsilon}

\def\Z{\ams{Z}}
\def\E{\ams{E}}
\def\H{\ams{H}}
\def\R{\ams{R}}
\def\C{\ams{C}}\def\Q{\ams{Q}}

\newfont{\Ma}{msam10}

\title{3-manifolds from Platonic solids}
\author{B. Everitt}
\address{Department of Mathematical Sciences, University of Durham, Durham DH1 3LE, England}
\email{brent.everitt@durham.ac.uk}
\curraddr{Department of Mathematics, University of York, York
YO10 5DD, England}
\email{bje1@york.ac.uk}


\thanks{The author would like to thank the Department of Mathematical
Sciences at the University of Aberdeen for the use of their
computational facilities.}

\begin{document}

\maketitle
\begin{abstract}
The problem of classifying, upto isometry (or similarity), the orientable spherical, Euclidean and
hyperbolic 3-manifolds that arise by identifying the faces of a Platonic solid
is formulated in the language of Coxeter groups.
In the spherical and hyperbolic cases, this allows us to complete the 
classification begun by Lorimer \cite{Lorimer92}, Richardson and Rubinstein
\cite{Rubinstein82} and Best \cite{Best71}.
\end{abstract}

\section{Introduction}\label{section_intro}

The first example of an orientable hyperbolic 3-manifold
arose by identifying the faces of a solid hyperbolic dodecahedron \cite{Weber33}.
Of course in the intervening years, much much more has been said about such manifolds.
Yet, the classical question of which spherical, Euclidean or hyperbolic manifolds arise by 
identifying the faces of a Platonic solid has
a surprisingly incomplete solution.

In this paper we formulate the problem in terms of
classifying certain subgroups of rank four Coxeter groups. 
This approach is implicit in \cite{Lorimer92,Rubinstein82}, and this
paper should be viewed as completing their work in the spherical and
hyperbolic cases. It follows an
earlier, oft quoted but flawed attempt in \cite{Best71}. As the results
of \cite{Rubinstein82} are not readily available in the literature, we
summarise them in Table \ref{table 2}.

There are other, non-algebraic, approaches to the problem, particularly that of Molna'r
and his school (see for instance \cite{Molnar92,Prok98,Prok96} and the
references there). In
fact, our list of manifolds in the Euclidean case cannot be given
precisely without recourse to Prok's paper \cite{Prok96}.
The author is very grateful to Colin Maclachlan, Emil Molna'r, Istvan
Prok, Peter Lorimer and Marston Conder for many useful
discussions and suggestions. He also thanks Hyam Rubinstein for a copy
of the preprint \cite{Rubinstein82}.

\section{Platonic solids and Coxeter
groups}\label{section_Coxeter_groups}

Let
$X=S^3,{\ams E}^3$ or
${\ams H}^3$, and suppose $\Delta\subset X$ is a finite volume Coxeter 
simplex (see \cite{Humphreys90}) with symbol, 
\begin{equation}\label{coxeter.symbol}
\begin{picture}(76,4)
\put(12,5){$p$}\put(36,5){$q$}\put(60,5){$r$}
\put(2,2){\circle{4}}\put(4,2){\line(1,0){20}}
\put(26,2){\circle{4}}\put(28,2){\line(1,0){20}}
\put(50,2){\circle{4}}\put(52,2){\line(1,0){20}}
\put(74,2){\circle{4}}\end{picture}
\end{equation}
Each node of the symbol corresponds to a face of $\Delta$, which in turn
has a
vertex of $\Delta$ opposite it. Call this the vertex corresponding to the node. 
Let $\Gamma=\{p,q,r\}$ be the Coxeter group generated by
reflections in the faces of $\Delta$, and for any vertex, edge or face of $\Delta$, say
$*$, let $\Gamma_*$ be its stabiliser in $\Gamma$. In particular, if $v$ is a vertex of $\Delta$, then
$\Gamma_v$ is also Coxeter group, its symbol obtained from (\ref{coxeter.symbol}) by deleting
the node corresponding to $v$ and its incident edges.

Let $v$ be the vertex of $\Delta$ corresponding to the left-most node of (\ref{coxeter.symbol}).
Then,
\begin{equation}\label{eq1}
\Sigma=\bigcup_{\gamma\in\Gamma_{v}} \gamma(\Delta),
\end{equation}
is a solid with $r$-gonal faces, $q$ meeting at each vertex, and dihedral
angle (that is, angle subtended by adjacent faces) $2\pi/p$.
Similarly for the last node with corresponding vertex $v'$, from which we obtain a solid
$\Sigma'$ with $p$-gonal faces, $q$
meeting at each vertex and dihedral angle $2\pi/r$. 
The tessellations of $X$ by congruent copies of
$\Sigma$ and $\Sigma'$ that result from
successive reflections in the faces of these solids are duals to one another, and both have
automorphism group $\Gamma$. 

On the otherhand, suppose we have a Platonic solid in $X$. By this we
mean a polytope $P$ with the combinatorial type of a Platonic solid
(convex regular solid),
embedded in $X$ so that all side lengths are equal, as are the interior
face angles and dihedral angles. For
face identifications of $P$ to yield an $X$-manifold, the dihedral angle
must be a submultiple of $2\pi$, say $2\pi/p$. Barycentric subdividion
of $P$ then gives a Coxeter simplex of the type (\ref{coxeter.symbol}),
and $P$ is recoverable in the form (\ref{eq1}) using the vertex $v$ of
the simplex lying at the center of $P$. Thus, the problem of obtaining
manifolds from a general Platonic solid $P$ reduces to consideration of the $\Sigma$
obtained at (\ref{eq1}).

All Coxeter simplices of the form (\ref{coxeter.symbol}) are known and listed in
Sections 2.4, 2.5 and 6.9 of
\cite{Humphreys90}. For
$X=S^3$ we have, 
$$
\begin{picture}(76,4)
\put(12,5){$$}\put(36,5){$$}\put(60,5){$$}
\put(2,2){\circle{4}}\put(4,2){\line(1,0){20}}
\put(26,2){\circle{4}}\put(28,2){\line(1,0){20}}
\put(50,2){\circle{4}}\put(52,2){\line(1,0){20}}
\put(74,2){\circle{4}}\end{picture},\,\,\,\,
\begin{picture}(76,4)
\put(12,5){$4$}\put(36,5){$$}\put(60,5){$$}
\put(2,2){\circle{4}}\put(4,2){\line(1,0){20}}
\put(26,2){\circle{4}}\put(28,2){\line(1,0){20}}
\put(50,2){\circle{4}}\put(52,2){\line(1,0){20}}
\put(74,2){\circle{4}}\end{picture},\,\,\,\,
\begin{picture}(76,4)
\put(12,5){$$}\put(36,5){$4$}\put(60,5){$$}
\put(2,2){\circle{4}}\put(4,2){\line(1,0){20}}
\put(26,2){\circle{4}}\put(28,2){\line(1,0){20}}
\put(50,2){\circle{4}}\put(52,2){\line(1,0){20}}
\put(74,2){\circle{4}}\end{picture},\,\,\,\,
\begin{picture}(76,4)
\put(12,5){$$}\put(36,5){$$}\put(60,5){$5$}
\put(2,2){\circle{4}}\put(4,2){\line(1,0){20}}
\put(26,2){\circle{4}}\put(28,2){\line(1,0){20}}
\put(50,2){\circle{4}}\put(52,2){\line(1,0){20}}
\put(74,2){\circle{4}}\end{picture};
$$
for $X={\ams E}^3$ we get
$$
\begin{picture}(76,4)
\put(12,5){$4$}\put(36,5){$$}\put(60,5){$4$}
\put(2,2){\circle{4}}\put(4,2){\line(1,0){20}}
\put(26,2){\circle{4}}\put(28,2){\line(1,0){20}}
\put(50,2){\circle{4}}\put(52,2){\line(1,0){20}}
\put(74,2){\circle{4}}\end{picture};
$$
and for
$X={\ams H}^3$,
\begin{align*}
&\begin{picture}(76,4)
\put(12,5){$4$}\put(36,5){$$}\put(60,5){$5$}
\put(2,2){\circle{4}}\put(4,2){\line(1,0){20}}
\put(26,2){\circle{4}}\put(28,2){\line(1,0){20}}
\put(50,2){\circle{4}}\put(52,2){\line(1,0){20}}
\put(74,2){\circle{4}}\end{picture},\,\,\,\,
\begin{picture}(76,4)
\put(12,5){$$}\put(36,5){$5$}\put(60,5){$$}
\put(2,2){\circle{4}}\put(4,2){\line(1,0){20}}
\put(26,2){\circle{4}}\put(28,2){\line(1,0){20}}
\put(50,2){\circle{4}}\put(52,2){\line(1,0){20}}
\put(74,2){\circle{4}}\end{picture},\,\,\,\,
\begin{picture}(76,4)
\put(12,5){$5$}\put(36,5){$$}\put(60,5){$5$}
\put(2,2){\circle{4}}\put(4,2){\line(1,0){20}}
\put(26,2){\circle{4}}\put(28,2){\line(1,0){20}}
\put(50,2){\circle{4}}\put(52,2){\line(1,0){20}}
\put(74,2){\circle{4}}\end{picture},\,\,\,\,
\begin{picture}(76,4)
\put(12,5){$4$}\put(36,5){$4$}\put(60,5){$$}
\put(2,2){\circle{4}}\put(4,2){\line(1,0){20}}
\put(26,2){\circle{4}}\put(28,2){\line(1,0){20}}
\put(50,2){\circle{4}}\put(52,2){\line(1,0){20}}
\put(74,2){\circle{4}}\end{picture},\\
\\
&\begin{picture}(76,4)
\put(12,5){$4$}\put(36,5){$$}\put(60,5){$6$}
\put(2,2){\circle{4}}\put(4,2){\line(1,0){20}}
\put(26,2){\circle{4}}\put(28,2){\line(1,0){20}}
\put(50,2){\circle{4}}\put(52,2){\line(1,0){20}}
\put(74,2){\circle{4}}\end{picture},\,\,\,\,
\begin{picture}(76,4)
\put(12,5){$5$}\put(36,5){$$}\put(60,5){$6$}
\put(2,2){\circle{4}}\put(4,2){\line(1,0){20}}
\put(26,2){\circle{4}}\put(28,2){\line(1,0){20}}
\put(50,2){\circle{4}}\put(52,2){\line(1,0){20}}
\put(74,2){\circle{4}}\end{picture},\,\,\,\,
\begin{picture}(76,4)
\put(12,5){$$}\put(36,5){$$}\put(60,5){$6$}
\put(2,2){\circle{4}}\put(4,2){\line(1,0){20}}
\put(26,2){\circle{4}}\put(28,2){\line(1,0){20}}
\put(50,2){\circle{4}}\put(52,2){\line(1,0){20}}
\put(74,2){\circle{4}}\end{picture},\,\,\,\,
\begin{picture}(76,4)
\put(12,5){$$}\put(36,5){$6$}\put(60,5){$$}
\put(2,2){\circle{4}}\put(4,2){\line(1,0){20}}
\put(26,2){\circle{4}}\put(28,2){\line(1,0){20}}
\put(50,2){\circle{4}}\put(52,2){\line(1,0){20}}
\put(74,2){\circle{4}}\end{picture},\\
\\
&\begin{picture}(76,4)
\put(12,5){$4$}\put(36,5){$4$}\put(60,5){$4$}
\put(2,2){\circle{4}}\put(4,2){\line(1,0){20}}
\put(26,2){\circle{4}}\put(28,2){\line(1,0){20}}
\put(50,2){\circle{4}}\put(52,2){\line(1,0){20}}
\put(74,2){\circle{4}}\end{picture},\,\,\,\,
\begin{picture}(76,4)
\put(12,5){$6$}\put(36,5){$$}\put(60,5){$6$}
\put(2,2){\circle{4}}\put(4,2){\line(1,0){20}}
\put(26,2){\circle{4}}\put(28,2){\line(1,0){20}}
\put(50,2){\circle{4}}\put(52,2){\line(1,0){20}}
\put(74,2){\circle{4}}\end{picture}\,\,\,\,
\end{align*}
In the spherical case, the tessellations of $S^3$ by copies of $\Sigma$
or $\Sigma'$ give the six 4-dimensional regular solids \cite{Hilbert52}.
In
another incarnation, the first three give $\Gamma$ that are the Weyl groups of the Lie algebras
of type $A_4=\s\l_{5}(\C)$, $B_4=\s\o_{9}(\C)$ and $F_4$. 
The hyperbolic $\Gamma$ give $\Sigma$ and $\Sigma'$ of finite volume: 
the first three compact, the others non-compact.

\begin{table}
\begin{center}\begin{tabular}{cccccc}
$N$&$FI$&$EI$&$H_1$\\
{\small{\tt 1}}&{\small{\tt abcdefefbcda}}&{\small{\tt
a(-+)b(-+)c(-+)d(-+)e(-+)f(-+)g(-+)}}&{\small{\tt 00000}}\\
&&{\small{\tt h(-+)i(-+)j(-+)idjefagbhcghijfeabcd}}&&\\

{\small{\tt 2}}&{\small{\tt abcdefbdcfea}}&{\small{\tt 
a(-+)b(-+)c(-+)d(-+)e(-+)f(++)g(++)}}&{\small{\tt (15)0000}}\\
&&{\small{\tt h(++)i(++)j(++)ajcgbfeidhfhgjieabcd}}\\ 
\end{tabular}\end{center}
\caption{The spherical manifolds arising from a dodecahedron
with dihedral angle $2\pi/3$,
\cite{Lorimer92}.}\label{table 1}\end{table}

We get a total of six spherical, one Euclidean and eight hyperbolic
Platonic solids from these groups: spherical tetrahedra with dihedral
angles $2\pi/3,2\pi/4$ and $2\pi/5$, a cube with angle $2\pi/3$, an
octahedron with angle $2\pi/3$ and a dodecahedron with angle $2\pi/3$;
in the Euclidean case we get the familiar cube; and in the hyperbolic,
a compact octahedron, icosahedron and two dodecahedrons with angles
$2\pi/5,2\pi/3,2\pi/4$ and $2\pi/5$; finally, a non-compact but finite
volume cube, octahedron, dodecahedron and tetrahedron with dihedral
angles $2\pi/4,2\pi/6,2\pi/6$ and $2\pi/6$ respectively.

\section{Constructing the manifolds}\label{sect3}

Any $X$-manifold (see \cite[\S 3.3]{Thurston97}) arises as the quotient $X/K$ of $X$ by a group $K$ acting
properly discontinuously and without fixed points. When $X=\E^3$ or $\H^3$, the isometries of
$X$ with fixed points are precisely those of finite order. This allows a
simple algebraic formulation of the problem in these two geometries
(Theorem \ref{thm1} below). Alternatively, 
recourse to a more geometric view yields Theorem \ref{thm2}, which holds
for all three geometries.
The statements in the remainder
of the paper will be formulated in terms of the solid $\Sigma$, those for
$\Sigma'$ being entirely analogous.

Establishing first some notation, let $\SS_m$ be the symmetric group of degree $m$.
If $\Lambda$ is a subgroup of $\SS_m$, let $\Lambda_i$ be the stabiliser
in $\Lambda$ of
$i\in\{1,\ldots,m\}$.
For any group $G$, let $\TT(G)$ be a subset that contains {\em at least one\/}
representative from each conjugacy class of elements of finite prime
order. 

\begin{theorem}\label{thm1}
Let $X=\E^n$ or $\H^n$ for $n\geq 2$; $\Gamma$ a group acting properly
discontinuously by isometries on $X$ with (convex, locally finite) fundamental region $P$; $F$ a finite
subgroup of $\Gamma$ and 
$$
\Sigma=\bigcup_{\gamma\in F} \gamma(P).
$$
An $X$-manifold $M$ arises by the identification of points on the
boundary of $\Sigma$
if and only if there is a homomorphism $\varepsilon:\Gamma\rightarrow
\SS_m$, where $m$ is the order of $F$, such that,
\begin{enumerate}
\item if $\Lambda=\varepsilon(\Gamma)$, then $\Lambda$ acts transitively
on $\{1,\ldots,m\}$, and
\item for all $\gamma\in\TT(\Gamma)$, the permutation $\varepsilon(\gamma)$
fixes no point of 
$\{1,\ldots,m\}$.
\end{enumerate}
Moreover, if $i\in\{1,\ldots,m\}$, 
then $\pi_1(M)\cong \varepsilon^{-1}(\Lambda_i)$.
\end{theorem}

\begin{proof}
An $X$-manifold $M$ arises by identifying points on $\partial\Sigma$ if and only if
there is a torsion free subgroup $K$ of $\Gamma$ 
with fundamental region $\Sigma$ and $M$ isometric to the
quotient $X/K$. 
Such a $K$ (which is isomorphic to $\pi_1(M)$) may be replaced by any of its
conjugates in $\Gamma$, as these will yield quotients isometric to $M$.
Conjugacy classes of subgroups of $\Gamma$ of index $m$ correspond to
transitive actions of $\Gamma$ on $\{1,\ldots,m\}$, the subgroups
arising as the stabilisers of points. These actions in turn correspond to
homomorphisms $\Gamma\rightarrow\SS_m$ with transitive image.
	

A subgroup $K$ is torsion free if and only if it intersects trivially
the conjugacy class in $\Gamma$ of each $\gamma\in\TT(\Gamma)$. This
happens precisely when $\ve(\gamma)$ has no fixed points among
$\{1,\ldots,m\}$. Finally, $\Sigma$ forms a fundamental region for the
action of $K$ on $X$ exactly when $F$ forms a transversal (a non-redundant
list of coset representatives) for $K$ in $\Gamma$. Equivalently, $K\cap
F=\{1\}$ and $KF=\Gamma$. The first follows immediatly as $K$ is torsion free, and
the second, since $F$ is a subgroup, when the index of $K$ in $\Gamma$
is equal to the order of $F$. 
\end{proof}

\begin{table}
\begin{center}\begin{tabular}{cccccc}
$N$&$FI$&$EI$&$H_1$\\
{\small{\tt 1}}&{\small{\tt abcdefefbcda}}&{\small{\tt
a(-+-+)b(-+-+)c(-+-+)d(-+-+)e(-+-+)}}&{\small{\tt 55500}}\\
&&{\small{\tt cdeabf(++++)afbfcfdfecdeabdeabc}}&&\\

{\small{\tt 2}}&{\small{\tt abcdefdefbca}}&{\small{\tt
a(++++)b(++++)c(++++)d(++++)e(++++)}}&{\small{\tt 55500}}\\
&&{\small{\tt abcdebf(++++)cfdfefafcdeabbcdea}}\\ 

{\small{\tt 3}}&{\small{\tt abcdefdefbca}}&{\small{\tt 
a(+-++)b(-+++)c(---+)d(++-+)e(+-++)}}&{\small{\tt 33000}}\\
&&{\small{\tt debaf(+-++)bcfafefcdcfedabeabcd}}\\ 

{\small{\tt 4}}&{\small{\tt abccadeefbfd}}&{\small{\tt
a(++--)ab(-+++)ac(-+-+)d(-+++)bab}}&{\small{\tt 57000}}\\
&&{\small{\tt e(+++-)ef(--+-)bfdcaecdfffddcbece}}\\
  
{\small{\tt 5}}&{\small{\tt abcdefebfdca}}&{\small{\tt
a(-+-+)b(-+-+)c(-+-+)d(-+-+)e(-+-+)}}&{\small{\tt 35500}}\\
&&{\small{\tt edacbf(++++)cfefbfafdbdaeceabcd}}\\
 
{\small{\tt 6}}&{\small{\tt abcdeffbdeca}}&{\small{\tt
a(++++)b(++++)c(++++)d(++++)e(++++)}}&{\small{\tt 33550}}\\
&&{\small{\tt cf(++++)efdfbfafeacdbdeabc}}\\ 

{\small{\tt 7}}&{\small{\tt abcdebedffca}}&{\small{\tt
a(+-++)b(+-++)c(---+)d(-+-+)e(-+++)}}&{\small{\tt 3(16)000}}\\
&&{\small{\tt cedaef(--+-)afdfbfcfebdcbacdeab}}\\ 

{\small{\tt 8}}&{\small{\tt abbcadefecfd}}&{\small{\tt
a(+++-)b(++-+)c(--++)ad(-++-)a}}&{\small{\tt (29)0000}}\\
&&{\small{\tt e(+-+-)dbbeadcf(+--+)acfceffdedbdbfc}}\\ 

{\small{\tt 9}}&{\small{\tt abcbdaefghihdefjgcji}}&{\small{\tt
a(-+)b(+-)c(--)d(-+)e(-+)deabf(++)}}&{\small{\tt (11)(11)000}}\\
&&{\small {\tt g(+-)h(-+)i(+-)iaccj(++)jhdebfgfghij}}\\

{\small{\tt 10}}&{\small{\tt abcdebfceghhiijjfgda}}&{\small{\tt
a(-+)b(-+)c(-+)d(--)e(++)cf(--)ea}}&{\small{\tt 90000}}\\
&&{\small{\tt g(--)ebh(+-)gi(++)dj(+-)fghhdiifjjabc}}\\

{\small{\tt 11}}&{\small{\tt abcdefbdgehiijjhfgca}}&{\small{\tt
a(++)b(++)c(++)d(++)e(+-)cdf(+-)ad}}&{\small{\tt 22900}}\\
&&{\small{\tt g(+-)bfh(-+)gi(+-)ej(-+)ijgjhehifabc}}\\

{\small{\tt 12}}&{\small{\tt abcdaefdgfhihcjjbige}}&{\small{\tt
a(++)b(+-)bc(+-)d(--)e(+-)baf(--)}}&{\small{\tt 57000}}\\
&&{\small{\tt g(+-)efgh(++)ghci(+-)dj(-+)jjdeiicahf}}\\

{\small{\tt 13}}&{\small{\tt abcdabefghcijidfjghe}}&{\small{\tt
a(++)ab(-+)c(++)d(++)e(--)bacf(+-)}}&{\small{\tt (29)0000}}\\
&&{\small{\tt g(+-)h(+-)ei(-+)j(++)djfidhgihebgjfc}}\\

{\small{\tt 14}}&{\small{\tt abcdaebdfghicjehjfgi}}&{\small{\tt
a(++)b(+-)bc(--)d(-+)e(++)bacdef(+-)}}&{\small{\tt (29)0000}}\\
&&{\small{\tt g(--)h(+-)di(-+)aj(--)ijfehgcighjf}}\\  
\end{tabular}\end{center}
\caption{The compact hyperbolic manifolds arising from a dodecahedron
with dihedral angle $2\pi/5$ and an icosahedron with angle $2\pi/3$,
\cite{Rubinstein82}.}\label{table 2}\end{table}

We will be applying Theorem \ref{thm1}
with $F$ the stabiliser $\Gamma_v$.
In an arbitrary Coxeter group $\Gamma$, a $\TT(\Gamma)$ can be
found using \cite{Bourbaki68,Howlett93}--list for example the conjugacy
class representatives in the maximal finite parabolic subgroups of $\Gamma$.
For the group with symbol
(\ref{coxeter.symbol}), or in fact for any 3-dimensional Euclidean or
hyperbolic Coxeter group, it is particularly easy to find a $\TT(\Gamma)$:
take the generating reflections and the powers of their pairwise
products that have prime order. 

More geometrically, suppose we have a subgroup $K$ of $\Gamma$ for which $\Gamma_v$ is a
transversal, and let $S$ be a face of $\Sigma$. In the tessellation of
$X$ by copies of $\Sigma$ there is a unique copy $\Sigma_S$ of $\Sigma$ with $\Sigma\cap\Sigma_S=S$. 
Since $\Sigma$
forms a fundamental region for $K$, there is a unique element $\gamma_S\in K$ sending
$\Sigma$ to $\Sigma_S$, and hence a face $S'$ of $\Sigma$ with $\gamma_S(S')=S$. The collection of
isometries $\{\gamma_S\}_{S\in\Sigma}$ yield a side-pairing of $\Sigma$ as in 
\cite[Section 10.1]{Ratcliffe94}. It follows
immediately from Theorems 10.1.2 and 10.1.3 of \cite{Ratcliffe94} that,

\begin{theorem}\label{thm2}
Let $X=S^3,\E^3$ or $\H^3$. An $X$-manifold $M$ arises by the
identification of faces of (\ref{eq1}) if and only if $\Gamma$ has a
subgroup $K$ of orientation preserving isometries, such that
\begin{enumerate}
\item $\Gamma_v$ forms a transversal in $\Gamma$ for $K$;
\item if $\{\gamma_S\}$ are the resulting side pairings of $\Sigma$,
then $\gamma_S$ fixes no point of $S'$; and
\item for $x\in\Sigma$, let $[x]$ denote the points of $\Sigma$ identified with it 
under the side pairing. If $x$ lies in the interior of an edge of $\Sigma$,  
then $[x]$ has cardinality $p$.
\end{enumerate}
\end{theorem}


So we merely require that the faces of $\Sigma$ are identified in pairs and the edges in groups of $p$.  
The identifications can be described algebraically as follows: since
$\Gamma$ acts transitively 
on the $k$-cells ($k=0,1,2,3$) of the tessellation of $X$ by $\Sigma$,
the faces of $\Sigma$ are in one to one correspondence with the cosets $(\Gamma_f)\gamma$,
where $f$ is the common face of $\Sigma$ and $\Delta$, and 
$\gamma\in\Gamma_v$. Two faces $(\Gamma_f)\gamma_1$ and
$(\Gamma_f)\gamma_2$ are identified by $K$ exactly when $(\Gamma_f)\gamma_1k=(\Gamma_f)\gamma_2$
for some
$k\in K$. Similarly for the edge identifications--take cosets of $\Gamma_e$ for $e$ the 
common edge
of $\Delta$ and $\Sigma$. 

We will say that two $X$-manifolds $M_1$ and $M_2$, for $X=S^3,\H^3$ (respectively $X=\E^3$)
are the same if and
only if there is an $X$-isometry (resp. $X$-similarity) between them. Equivalently, if $M_i=X/K_i$,
then the $K_i$ are conjugate in the group of isometries (resp.
similarities) of $X$. In the hyperbolic case, the following will help in
distinguishing manifolds:

\begin{theorem}[Margulis \cite{Borel81}]\label{margulis}
Let $G$ be a connected semisimple Lie group with trivial centre and no
compact factors, and $\Gamma$ an irreducible lattice (discrete subgroup
with finite Haar measure) in $G$. Then the
commensurator $\text{comm}(\Gamma)$ is discrete if and only if $\Gamma$
is non-arithmetic.
\end{theorem}

Arithmetic is meant here in the sense of \cite{Borel62}, and the
commensurator of $\Gamma$ is 
the subgroup consisting of those $h\in G$ such that $\Gamma$ and
$h^{-1}\Gamma h$ are commensurable (have intersection of finite index in
each). 
If we take $G=PO_{1,3}(\R)$ to be the full isometry group of $\H^3$, then $G$ has two
connected components, one of which, $G^+=PSO_{1,3}(\R)\cong PSL_2(\C)$, consists of the orientation
preserving isometries. Although $G$ is thus not connected, it is nevertheless easy to see that Theorem \ref{margulis}
holds for $\Gamma$ in $G$. The arithmeticity of hyperbolic Coxeter
groups is easily determined using the results of \cite{Vinberg67}, from
which we get in particular that the group with symbol
\begin{equation}\label{536}
\begin{picture}(76,6)
\put(12,5){$5$}\put(36,5){$$}\put(60,5){$6$}
\put(2,2){\circle{4}}\put(4,2){\line(1,0){20}}
\put(26,2){\circle{4}}\put(28,2){\line(1,0){20}}
\put(50,2){\circle{4}}\put(52,2){\line(1,0){20}}
\put(74,2){\circle{4}}\end{picture},
\end{equation}
is {\em non}-arithmetic.
In \cite{Adams92} the six cofinite discrete subgroups of $G$
having the smallest covolume are enumerated. In particular, they
are all commensurable with the Bianchi groups $PGL_2{\mathcal O}_1$ or 
$PGL_2{\mathcal O}_3$, where
${\mathcal O}_d$ is the ring of integers in the number field
$\Q(\sqrt{-d})$. By comparing volumes, one sees that if the group
$\Gamma$ with
symbol (\ref{536}) is not maximal, then it is contained in a $\overline{\Gamma}$
one of the six above, which cannot be, for the six above are
arithmetic, while $\overline{\Gamma}$, being commensurable with $\Gamma$, is not.


Suppose then we have $K_i,i=1,2$, torsion free subgroups of the $\Gamma$
with symbol (\ref{536}), and
a $g\in PO_{1,3}(\R)$ such that $g^{-1}K_1g=K_2$. Then $g\in\text{comm}(\Gamma)$. By Theorem \ref{margulis},
$\text{comm}(\Gamma)$ is also discrete in $G$, and by the maximality of
$\Gamma$, we
have $\Gamma=\text{comm}(\Gamma)$. Thus $g\in\Gamma$. This reduces
consideration of the conjugacy of the $K_i$ in $G$ (which is
hard), to the much easier question of their conjugacy in $\Gamma$.

\begin{figure}\begin{center}
\begin{tabular}{ccccc}
\begin{pspicture}(0,0)(3,3)
\rput(2.2,2.2){${\scriptstyle 1}$}\rput(.8,2.2){${\scriptstyle 2}$}
\rput(1.65,.85){${\scriptstyle 3}$}
\rput(1.5,2.8){${\scriptstyle 4}$}
\rput(.6,1.3){${\scriptstyle 5}$}\rput(2.4,1.3){${\scriptstyle 6}$}
\rput(1.5,2){{\bf 1}}\rput(1.1,1.4){{\bf 2}}\rput(1.9,1.4){{\bf 3}}
\rput(1.5,3.2){{\bf 4}}
\psline{->}(0,2.59)(1.6,2.59)\psline(1.5,2.59)(3,2.59)
\psline{->}(1.5,0)(.675,1.43)\psline{->}(3,2.6)(2.175,1.17)
\psline{->}(1.5,1.5)(.675,2.105)\psline{->}(1.5,1.5)(2.325,2.105)
\psline{->}(1.5,1.5)(1.5,.75)
\psline(2.25,2.05)(3,2.6)\psline(1.5,0.75)(1.5,0)\psline(2.25,1.3)(1.5,0)
\psline(.75,1.3)(0,2.6)\psline(.75,2.05)(0,2.6)
\end{pspicture}&\hspace*{0.5cm}&
\begin{pspicture}(0,0)(4,5)
\rput(2,2){{\bf 1}}\rput(2,3.5){{\bf 2}}\rput(3.5,2){{\bf 3}}
\rput(2,0.5){{\bf 4}}\rput(0.5,2){{\bf 5}}\rput(2,4.8){{\bf 6}}


\rput(2.8,2){${\scriptstyle 1}$}\rput(2,2.8){${\scriptstyle 2}$}
\rput(3.5,3.25){${\scriptstyle 3}$}\rput(2,1.2){${\scriptstyle 4}$}
\rput(3.55,0.7){${\scriptstyle 5}$}\rput(4.2,2){${\scriptstyle 6}$}
\rput(1.2,2){${\scriptstyle 7}$}\rput(0.5,3.25){${\scriptstyle 8}$}
\rput(0.45,0.7){${\scriptstyle 9}$}\rput(2,4.2){${\scriptstyle 10}$}
\rput(2,-0.2){${\scriptstyle 11}$}\rput(-0.3,2){${\scriptstyle 12}$}

\psline{->}(0,4)(2.1,4)\psline(2,4)(4,4)
\psline{->}(1,3)(2.1,3)\psline(2,3)(3,3)
\psline{->}(1,1)(2.1,1)\psline(2,1)(3,1)
\psline{->}(0,0)(2.1,0)\psline(2,0)(4,0)

\psline{->}(0,0)(0,2.1)\psline(0,2)(0,4)
\psline{->}(1,1)(1,2.1)\psline(1,2)(1,3)
\psline{->}(3,1)(3,2.1)\psline(3,2)(3,3)
\psline{->}(4,0)(4,2.1)\psline(4,2)(4,4)

\psline{->}(1,3)(0.4,3.6)\psline(0.5,3.5)(0,4)
\psline{->}(3,3)(3.6,3.6)\psline(3.5,3.5)(4,4)
\psline{->}(3,1)(3.6,0.4)\psline(3.5,0.5)(4,0)
\psline{->}(1,1)(0.4,0.4)\psline(0.5,0.5)(0,0)
\end{pspicture}&\hspace*{0.5cm}&
\begin{pspicture}(0,0)(5,5)
\rput(2.5,2.3){${\scriptstyle 1}$}\rput(2,3.1){${\scriptstyle 2}$}
\rput(3,3.1){${\scriptstyle 9}$}\rput(1,3.55){${\scriptstyle 3}$}
\rput(4,3.55){${\scriptstyle 5}$}
\rput(2.3,1.3){${\scriptstyle 4}$}\rput(2.65,1.3){${\scriptstyle 12}$}
\rput(1.2,3.8){${\scriptstyle 10}$}\rput(3.8,3.8){${\scriptstyle 6}$}
\rput(2.5,4.5){${\scriptstyle 7}$}\rput(3.85,2){${\scriptstyle 8}$}
\rput(1.1,2){${\scriptstyle 11}$}

\rput(2.5,2.9){{\bf 1}}\rput(1.6,3.4){{\bf 2}}\rput(3.4,3.4){{\bf 5}}
\rput(2.5,1.8){{\bf 4}}\rput(2.5,4){{\bf 6}}\rput(2.5,5){{\bf 7}}
\rput(1.5,2.4){{\bf 3}}\rput(3.5,2.4){{\bf 8}}

\psline{->}(2.5,0)(1.125,2.3815)
\psline{->}(0,4.33)(2.75,4.33)
\psline{->}(5,4.33)(3.625,1.9485)
\psline{->}(1.8,2.5)(.81,3.5065)\psline{->}(1.8,2.5)(2.185,1.125)
\psline{->}(1.8,2.5)(2.185,3.1325)
\psline{->}(2.5,3.65)(1.125,4.024)\psline{->}(2.5,3.65)(3.875,4.024)
\psline{->}(2.5,3.65)(2.885,3.0175)
\psline{->}(3.2,2.5)(2.43,2.5)\psline{->}(3.2,2.5)(4.19,3.5065)
\psline{->}(3.2,2.5)(2.815,1.125)
\psline(2.15,1.25)(2.5,0)\psline(2.85,1.25)(2.5,0)\psline(3.75,2.165)(2.5,0)
\psline(.9,3.415)(0,4.33)\psline(1.25,2.165)(0,4.33)\psline(1.25,3.99)(0,4.33)
\psline(2.5,4.33)(5,4.33)\psline(3.75,3.99)(5,4.33)\psline(4.1,3.415)(5,4.33)
\psline(2.5,2.5)(1.8,2.5)\psline(2.15,3.075)(2.5,3.65)\psline(2.85,3.075)(3.2,2.5)
\end{pspicture}\\
\vspace*{1em}&&&&\\
\end{tabular}
\end{center}\caption{}\label{fig1}
\end{figure}



\section{The Manifolds}\label{section_spherical}

Of the fifteen Platonic solids listed at the end of Section \ref{section_intro}, four
can be removed from consideration using Theorem \ref{thm2}, as the number of
edges of $\Sigma$ is not divisible by $p$. Of those that remain, the
spherical dodecahedron with dihedral angle $2\pi/3$ was handled in \cite{Lorimer92}
with results listed in Table \ref{table 1} (the notation is described
below). The first of the two manifolds is the Poincar\'{e} homology sphere.
The compact hyperbolic dodecahedron and icosahedron with angles $2\pi/5$
and $2\pi/3$ were investigated in \cite{Rubinstein82} with the results
in Table \ref{table 2}--the first eight manifolds come from the
dodecahedron, the remainder from the icosahedron\footnote{It should be noted that while there are pairs in Table
\ref{table 2} with the same first homology, algebraic arguments are
provided in \cite{Rubinstein82} that show that the list is non-redundant
(this is to be contrasted with the list in \cite{Best71} which contains
isometric pairs). Generally this involves consideration of quotients
of terms in the derived series for $K=\pi_1(M)$, for instance, $K'/K''$.}. 
The first is the Weber-Seifert 
space.
This leaves the spherical $\{3,3,3\},\{4,3,3\}$ and $\{3,4,3\}$; the
Euclidean $\{4,3,4\}$ and hyperbolic $\{4,4,3\},\{4,3,6\},\{5,3,6\}$
and $\{3,3,6\}$.

As no doubt the reader has gathered by now, the only practical way the techniques
of the previous
section can be
implemented is computationally. 
We use Sims's low index subgroups algorithm as implemented in
Magma \cite{Cannon84} to find the homomorphisms required by Theorem \ref{thm1} when $X=\E^3$
and $\H^3$. For the spherical manifolds, we use Theorem \ref{thm2}.
In any case, we obtain  a complete list of the $K$, subgroups of 
the various $\Gamma$, satisfying the conditions of the two Theorems. As we
want orientable manifolds, we also require that the generators of $K$
are words of even length in the generators for $\Gamma$. 
The resulting $K$ will be {\em non-conjugate in $\Gamma$\/}, although not
necessarily so in $G$, the full isometry/similarity group of $X$. 


\begin{figure}\begin{center}
\begin{tabular}{ccc}
\begin{pspicture}(5,6)

\rput(2.5,2.9){{\bf 1}}\rput(3.15,3.75){{\bf 2}}
\rput(1.85,3.75){{\bf 6}}\rput(3.65,2.6){{\bf 3}}
\rput(1.35,2.6){{\bf 5}}\rput(4.4,3.4){{\bf 7}}
\rput(.6,3.4){{\bf 10}}\rput(2.5,1.8){{\bf 4}}
\rput(2.5,4.8){{\bf 11}}\rput(2.5,6){{\bf 12}}
\rput(3.7,1.4){{\bf 8}}\rput(1.3,1.4){{\bf 9}}

\rput(2.7,3.3){${\scriptstyle 1}$}\rput(2.3,3.3){${\scriptstyle 5}$}
\rput(2.95,2.8){${\scriptstyle 2}$}\rput(2.05,2.8){${\scriptstyle 4}$}
\rput(2.5,2.5){${\scriptstyle 3}$}
\rput(3.5,3.35){${\scriptstyle 6}$}\rput(1.5,3.35){${\scriptstyle 9}$}
\rput(3.25,2.25){${\scriptstyle 7}$}\rput(1.75,2.25){${\scriptstyle 8}$}
\rput(2.7,3.9){${\scriptstyle 10}$}
\rput(4,3.9){${\scriptstyle 11}$}\rput(1,3.9){${\scriptstyle 18}$}
\rput(4.275,2.9){${\scriptstyle 12}$}\rput(.725,2.9){${\scriptstyle
17}$}
\rput(3.9,1.95){${\scriptstyle 13}$}\rput(1.1,1.95){${\scriptstyle 16}$}
\rput(3.1,1.4){${\scriptstyle 14}$}\rput(1.9,1.4){${\scriptstyle 15}$}
\rput(1.9,4.5){${\scriptstyle 19}$}\rput(3.1,4.5){${\scriptstyle 20}$}
\rput(3.8,4.9){${\scriptstyle 21}$}\rput(1.2,4.9){${\scriptstyle 25}$}
\rput(4.8,2.4){${\scriptstyle 22}$}\rput(.2,2.4){${\scriptstyle 24}$}
\rput(2.7,.6){${\scriptstyle 23}$}
\rput(5.1,3.6){${\scriptstyle 27}$}\rput(-.1,3.6){${\scriptstyle 30}$}
\rput(4.15,.9){${\scriptstyle 28}$}\rput(0.85,.9){${\scriptstyle 29}$}
\rput(2.5,5.6){${\scriptstyle 26}$}

\psline{->}(5.3,1.95)(3.76,.8775)
\psline{->}(2.5,1.1)(2.5,.495)
\psline{->}(3.35,1.875)(2.8825,1.44875)
\psline{->}(2.925,2.35)(3.15875,2.08875)
\psline{->}(2.5,1.1)(2.0325,1.52625)
\psline{->}(4.3,2.35)(4.85,2.13)
\psline{->}(4.3,2.35)(3.7775,2.08875)
\psline{->}(2.5,0)(0.96,1.0725)
\psline{->}(2.075,2.35)(1.84125,2.08875)
\psline{->}(1.65,1.875)(1.1275,2.13625)
\psline{->}(-.3,1.95)(0.2225,3.8475)
\psline{->}(.7,2.35)(.15,2.13)
\psline{->}(.7,2.35)(.93375,2.8725)
\psline{->}(2.075,2.35)(1.92375,2.7625)
\psline{->}(2.925,2.35)(2.45,2.35)
\psline{->}(3.2,3.1)(3.04875,2.6875)
\psline{->}(3.2,3.1)(3.57125,3.21)
\psline{->}(1.8,3.1)(2.185,3.4025)
\psline{->}(1.8,3.1)(1.42875,3.21)
\psline{->}(1.125,3.3)(1.22125,3.905)
\psline{->}(3.875,3.3)(4.10875,2.7775)
\psline{->}(2.5,3.65)(2.885,3.3475)
\psline{->}(2.5,3.65)(2.5,4.0075)
\psline{->}(1.3,4.4)(1.96,4.345)
\psline{->}(1.3,4.4)(.9425,4.95)
\psline{->}(2.5,4.3)(3.16,4.355)
\psline{->}(3.7,4.4)(3.79625,3.795)
\psline{->}(3.7,4.4)(4.0575,4.95)
\psline{->}(.65,5.4)(2.685,5.4)
\psline{->}(4.35,5.4)(4.8725,3.5025)

\psline(2.075,1.4875)(1.65,1.875)
\psline(1.8625,2.1125)(1.65,1.875)
\psline(2.925,1.4875)(2.5,1.1)
\psline(3.825,2.1125)(3.35,1.875)
\psline(3.1375,2.1125)(3.35,1.875)
\psline(2.5,.55)(2.5,0)
\psline(3.9,.975)(2.5,0)
\psline(1.1,.975)(-.3,1.95)
\psline(.2,2.15)(-.3,1.95)
\psline(4.8,2.15)(5.3,1.95)
\psline(4.825,3.675)(5.3,1.95)
\psline(1.175,2.1125)(.7,2.35)
\psline(4.0875,2.825)(4.3,2.35)
\psline(2.5,2.35)(2.075,2.35)
\psline(3.0625,2.725)(2.925,2.35)
\psline(1.4625,3.2)(1.125,3.3)
\psline(.9125,2.825)(1.125,3.3)
\psline(1.9375,2.725)(1.8,3.1)
\psline(2.85,3.375)(3.2,3.1)
\psline(3.5375,3.2)(3.875,3.3)
\psline(3.7875,3.85)(3.875,3.3)
\psline(2.15,3.375)(2.5,3.65)
\psline(1.2125,3.85)(1.3,4.4)
\psline(1.9,4.35)(2.5,4.3)
\psline(2.5,3.975)(2.5,4.3)
\psline(3.1,4.35)(3.7,4.4)
\psline(.175,3.675)(.65,5.4)
\psline(.975,4.9)(.65,5.4)
\psline(2.5,5.4)(4.35,5.4)
\psline(4.825,3.675)(4.35,5.4)
\psline(4.025,4.9)(4.35,5.4)
\end{pspicture}&\hspace*{1cm}&
\begin{pspicture}(5,6)
\rput(2.5,3.4){{\bf 1}}\rput(3.4,3.65){{\bf 2}}
\rput(3.85,3){{\bf 3}}\rput(3.5,2.35){{\bf 4}}\rput(2.5,2.3){{\bf 5}}
\rput(1.5,2.35){{\bf 6}}\rput(1.15,3){{\bf 7}}\rput(1.6,3.65){{\bf 8}}
\rput(2,4.45){{\bf 9}}\rput(3,4.45){{\bf 10}}
\rput(4.4,4.65){{\bf 11}}\rput(4.65,4.25){{\bf 12}}
\rput(2.75,.9){{\bf 13}}\rput(2.25,.9){{\bf 14}}
\rput(.325,4.25){{\bf 15}}\rput(.6,4.65){{\bf 16}}
\rput(2.5,5.325){{\bf 17}}\rput(2.5,6.1){{\bf 20}}
\rput(4.275,2.3){{\bf 18}}\rput(.725,2.3){{\bf 19}}

\rput(2.8,3.4){${\scriptstyle 1}$}\rput(2.5,2.95){${\scriptstyle 2}$}
\rput(2.175,3.4){${\scriptstyle 3}$}\rput(3.55,4.05){${\scriptstyle 4}$}
\rput(3.75,3.6){${\scriptstyle 5}$}\rput(3.725,2.75){${\scriptstyle 6}$}
\rput(2.8,2){${\scriptstyle 7}$}\rput(2.2,2){${\scriptstyle 8}$}
\rput(1.275,2.75){${\scriptstyle 9}$}\rput(1.275,3.6){${\scriptstyle 10}$}
\rput(1.45,4.05){${\scriptstyle 11}$}\rput(2.325,4.45){${\scriptstyle 12}$}
\rput(3.65,4.75){${\scriptstyle 13}$}\rput(4.025,3.4){${\scriptstyle 14}$}
\rput(3.1,1.85){${\scriptstyle 15}$}\rput(1.325,4.75){${\scriptstyle 18}$}
\rput(.95,3.4){${\scriptstyle 17}$}
\rput(1.9,1.85){${\scriptstyle 16}$}\rput(3.85,5.05){${\scriptstyle 19}$}
\rput(4.725,4.725){${\scriptstyle 20}$}\rput(4.85,3.65){${\scriptstyle 21}$}
\rput(3.45,1.15){${\scriptstyle 22}$}\rput(2.675,.62){${\scriptstyle 23}$}
\rput(1.15,5.05){${\scriptstyle 27}$}
\rput(.275,4.725){${\scriptstyle 26}$}\rput(.175,3.65){${\scriptstyle 25}$}
\rput(1.55,1.15){${\scriptstyle 24}$}
\rput(2.5,5.75){${\scriptstyle 28}$}
\rput(4.75,2.1){${\scriptstyle 29}$}\rput(.25,2.1){${\scriptstyle 30}$}

\psline{->}(2.5,1.1)(1.675,1.815)\psline{->}(2.5,1.1)(2.5,.495)
\psline{->}(4,2.4)(3.175,1.08)\psline{->}(4,2.4)(3.175,1.685)
   \psline{->}(4,2.4)(4.6875,3.7475)
\psline{->}(3.4,2.8)(2.41,2.8)\psline{->}(3.4,2.8)(2.905,1.865)
   \psline{->}(3.4,2.8)(3.9775,3.6525)\psline{->}(3.4,2.8)(3.73,2.58)
\psline{->}(4.45,4.35)(4.2025,3.2775)\psline{->}(4.45,4.35)(4.89,4.625)
\psline{->}(2.5,4.05)(1.4275,4.215)\psline{->}(2.5,4.05)(2.5,4.545)
  \psline{->}(2.5,4.05)(3.5725,4.215)\psline{->}(2.5,4.05)(2.995,3.3625)
\psline{->}(2.5,4.95)(.9875,4.895)\psline{->}(2.5,4.95)(4.0125,4.895)
  \psline{->}(2.5,4.95)(3.5725,4.62)
\psline{->}(.55,4.35)(.11,4.625)\psline{->}(.55,4.35)(1.6225,4.68)
\psline{->}(1.6,2.8)(1.0225,3.6525)\psline{->}(1.6,2.8)(1.27,2.58)
  \psline{->}(1.6,2.8)(2.095,1.865)\psline{->}(1.6,2.8)(2.095,3.4875)
\psline{->}(1,2.4)(.3125,3.7475)\psline{->}(1,2.4)(.7525,3.4725)
  \psline{->}(1,2.4)(1.825,1.08)

\psline{-}(1.75,1.2)(2.5,0)\psline{-}(2.5,.55)(2.5,0)
  \psline{-}(3.25,1.2)(2.5,0)
\psline{-}(2.05,1.95)(2.5,1.1)\psline{-}(2.95,1.95)(2.5,1.1)
  \psline{-}(3.25,1.75)(2.5,1.1)
\psline{-}(3.7,2.6)(4,2.4)\psline{-}(4.225,3.375)(4,2.4)
\psline{-}(2.95,3.425)(3.4,2.8)
\psline{-}(3.475,4.65)(4.45,4.35)\psline{-}(3.475,4.2)(4.45,4.35)
  \psline{-}(3.925,3.575)(4.45,4.35)
\psline{-}(3.875,4.9)(5.25,4.85)\psline{-}(4.85,4.6)(5.25,4.85)
  \psline{-}(4.625,3.625)(5.25,4.85)
\psline{-}(2.05,3.425)(2.5,4.05)
\psline{-}(1.525,4.65)(2.5,4.95)\psline{-}(2.5,4.5)(2.5,4.95)
\psline{-}(.375,3.625)(-.25,4.85)\psline{-}(.15,4.6)(-.25,4.85)
  \psline{-}(1.125,4.9)(-.25,4.85)
\psline{-}(.775,3.375)(.55,4.35)\psline{-}(1.075,3.575)(.55,4.35)
  \psline{-}(1.525,4.2)(.55,4.35)
\psline{-}(2.5,2.8)(1.6,2.8)
\psline{-}(1.3,2.6)(1,2.4)\psline{-}(1.75,1.75)(1,2.4)

\pscurve{-}(2.5,0)(4.55,2.1)(5.25,4.85)\psline{->}(4.55,2.1)(4.5,2)
\pscurve{-}(5.25,4.85)(2.5,5.6)(-.25,4.85)\psline{->}(2.5,5.6)(2.53,5.6)
\pscurve{-}(-.25,4.85)(.45,2.1)(2.5,0)\psline{->}(.45,2.1)(.4,2.2)
\end{pspicture}\\
\end{tabular}
\caption{}\label{fig2}
\end{center}\end{figure}

The results are listed in Tables \ref{table 3}-\ref{table 5} which we will discuss in some detail presently.
First we describe the notation. In each of the Tables, the column headed $N$ indexes the manifolds
$M_i$ carrying the indicated geometric structure. The columns $FI$ and $EI$ give the face and edge
identifications in the form of an encoded string of letters and $\pm $
signs to be read in conjunction with Figures
\ref{fig1}-\ref{fig2}. The $i$-th and $j$-th faces are paired when the $i$-th and $j$-th positions of the
string in column $FI$ are occupied by the same letter. 
Similarly for the edge identifications, where a string of $\pm$'s after
a letter indicates whether the corresponding edge is identified with
subsequent ones with the orientations matching or reversed. For example,
the manifold $M_{18}$ arising from the dodecahedron $\{5,3,6\}$ has edge identifications
$$
\begin{tabular}{c}
{\tt a(+---+)b(+--++)bc(++--++)d(---+-)}\\
{\tt bcae(+-+--)ceadddbeacedcaabecbed},\\
\end{tabular}
$$
where {\tt e} indicates that edges $9,11,17,20,26$ and $29$ are identified,
and {\tt e(+-+--)} means edge $9$ is identified with edge $11$ so
that the identifications match, with edge $17$ so they are reversed,
with edge $20$ so they match, and so on. From the data in these two
columns one may reconstruct the side pairing transformations
$\{\gamma_S\}_{s\in\Sigma}$.
In particular, the vertex identifications can be obtained in the
spherical and Euclidean cases; in the hyperbolic there are no vertices!
(they lie on the boundary of hyperbolic space in these non-compact examples).
The next column in Table \ref{table 5}
gives the number of cusps. The
final column gives the first homology
$H_1(M_i,\Z)=\Z_a\oplus\Z_b\oplus\Z_c\oplus\Z_d\oplus\Z^e$ in the form of a sequence
{\tt abcde} (brackets are used in Tables \ref{table 1}-\ref{table 2} to
distinguish double digits).

\begin{table}
\begin{center}\begin{tabular}{cccc}
$N$&$FI$&$EI$&$H_1$\\
{\small{\tt 1}}&{\small{\tt abab}}&{\small{\tt a(--)b(--)aabb}}&{\small{\tt 50000}}\\
{\small{\tt 2}}&{\small{\tt ababcc}}&{\small{\tt a(++)b(+-)aac(+-)bcd(+-)bcdd}}&{\small{\tt 80000}}\\
{\small{\tt 3}}&{\small{\tt abcbca}}&{\small{\tt a(++)b(--)c(+-)cd(--)bdabdac}}&{\small{\tt 22000}}\\
{\small{\tt 4}}&{\small{\tt abcacbdd}}&{\small{\tt a(++)b(+-)c(+-)ad(++)cbdacdb}}&{\small{\tt 26000}}\\
{\small{\tt 5}}&{\small{\tt abcacdbd}}&{\small{\tt a(++)b(-+)c(++)ad(-+)cbcaddb}}&{\small{\tt 80000}}\\
{\small{\tt 6}}&{\small{\tt abcdcdab}}&{\small{\tt a(++)b(++)c(++)d(++)bcdadabc}}&{\small{\tt 30000}}\\
\end{tabular}\end{center}
\caption{The spherical manifolds}\label{table 3}\end{table}

\begin{table}
\begin{center}\begin{tabular}{cccc}
$N$&$FI$&$EI$&$H_1$\\
{\small{\tt 7}}&{\small{\tt abacbc}}&{\small{\tt a(+++)b(+++)aac(+++)bccbcba}}&{\small{\tt 30001}}\\
{\small{\tt 8}}&{\small{\tt abbcca}}&{\small{\tt a(-+-)ab(--+)c(-+-)bacbbacc}}&{\small{\tt 22001}}\\
{\small{\tt 9}}&{\small{\tt abccba}}&{\small{\tt a(-+-)ab(--+)c(+--)bccbbcaa}}&{\small{\tt 44000}}\\
{\small{\tt 10}}&{\small{\tt abcbca}}&{\small{\tt a(+++)b(+++)c(+++)bcaaccbba}}&{\small{\tt 00003}}\\
{\small{\tt 11}}&{\small{\tt abcbca}}&{\small{\tt a(+++)b(+++)c(-+-)cbaacbbca}}&{\small{\tt 20001}}\\
{\small{\tt 12}}&{\small{\tt abcbca}}&{\small{\tt a(-+-)b(+--)c(+++)bcaaccbba}}&{\small{\tt 22001}}\\
\end{tabular}\end{center}
\caption{The Euclidean manifolds}\label{table 4}\end{table}

\begin{table}
\begin{center}\begin{tabular}{cccccc}
$N$&$FI$&$EI$&$C$&$H_1$\\
{\small{\tt 13}}&{\small{\tt ababcdcd}}&{\small{\tt a(---)aaab(---)c(+-+)bccbcb}}&{\small{\tt 2}}&{\small{\tt
00002}}\\ 
{\small{\tt 14}}&{\small{\tt abacbdcd}}&{\small{\tt a(-+-)b(+--)babbaac(---)ccc}}&{\small{\tt 2}}&{\small{\tt
00002}}\\ 
{\small{\tt 15}}&{\small{\tt ababcc}}&{\small{\tt a(++---)b(-+-+-)aabbbbaaba}}&{\small{\tt 2}}&{\small{\tt
20002}}\\ 
{\small{\tt 16}}&{\small{\tt ababcc}}&{\small{\tt a(++--+)b(+++--)aabbbabbaa}}&{\small{\tt 1}}&{\small{\tt
24001}}\\ 
{\small{\tt 17}}&{\small{\tt abcbca}}&{\small{\tt a(+-+-+)b(----+)bbabaabaab}}&{\small{\tt 2}}&{\small{\tt
20002}}\\ 
{\small{\tt 18}}&{\small{\tt abacbddceeff}}&{\small{\tt
a(+---+)b(+--++)bc(++--++)d(---+-)}}&{\small{\tt 1}}&{\small{\tt 20001}}\\
&&{\small{\tt bcae(+-+--)ceadddbeacedcaabecbed}}&&\\
{\small{\tt 19}}&{\small{\tt abacdcdbefef}}&{\small{\tt 
a(--+-+)b(+-+--)bc(+----)d(--+-+)}}&{\small{\tt 2}}&{\small{\tt 20002}}\\
&&{\small{\tt bcacdcadddae(----+)badcaceeeebeb}}\\ 
{\small{\tt 20}}&{\small{\tt abacdbdcefef}}&{\small{\tt 
a(-+-++)b(+--++)bc(+--+-)d(+---+)}}&{\small{\tt 2}}&{\small{\tt 20002}}\\
&&{\small{\tt bcaadcadddce(----+)bcdacaeeeebeb}}\\ 
{\small{\tt 21}}&{\small{\tt abcacdedeffb}}&{\small{\tt
a(++---)b(--+++)abbc(+-++-)bbc}}&{\small{\tt 1}}&{\small{\tt 22001}}\\
&&{\small{\tt d(--++-)dadbadde(---+-)ceecceeadedc}}\\  
{\small{\tt 22}}&{\small{\tt abcacdedfebf}}&{\small{\tt
a(+--++)b(--+++)abbc(+-+++)bbc}}&{\small{\tt 1}}&{\small{\tt 22002}}\\
&&{\small{\tt d(-+++-)e(+-+--)edbadedeadcceeaedac}}\\ 
{\small{\tt 23}}&{\small{\tt abcacdedfefb}}&{\small{\tt
a(+-++-)b(--+++)abbc(++++-)bbc}}&{\small{\tt 1}}&{\small{\tt 26002}}\\
&&{\small{\tt d(-+-+-)e(+-+++)edbadecdacdceeaecad}}\\ 
{\small{\tt 24}}&{\small{\tt abcacdedeffb}}&{\small{\tt
a(+----)b(--+++)abbc(++-+-)bbc}}&{\small{\tt 2}}&{\small{\tt 22002}}\\
&&{\small{\tt d(-++-+)ae(---+-)dbadecdecdceeeacad}}\\ 
{\small{\tt 25}}&{\small{\tt abcacbdeedff}}&{\small{\tt
a(+-+--)b(+-+++)ac(+--++)d(-++-+)}}&{\small{\tt 2}}&{\small{\tt 60002}}\\
&&{\small{\tt e(--+-+)dbdedbccaebeadceecbabacd}}\\ 
{\small{\tt 26}}&{\small{\tt abcacdebdeff}}&{\small{\tt
a(+--++)b(+-+++)ac(-+-+-)d(-++--)}}&{\small{\tt 2}}&{\small{\tt 20002}}\\
&&{\small {\tt e(+-+--)dbdecbccaebeacdeedbabadc}}\\
{\small{\tt 27}}&{\small{\tt abcbdefdcfae}}&{\small{\tt
a(-++--)b(--+++)ac(+++--)d(--+--)}}&{\small{\tt 1}}&{\small{\tt 22001}}\\
&&{\small{\tt ccbdbdae(+++-+)daeeebaceccadebbd}}\\  
\end{tabular}\end{center}
\caption{The hyperbolic manifolds}\label{table 5}\end{table}

Table 1 gives the spherical results. Manifold $M_1$ comes from the tetrahedron in $\{3,3,3\}$,
$M_2$ and $M_3$ from the cube in $\{4,3,3\}$ and $M_4,M_5$ and $M_6$ from the octahedron in
$\{3,4,3\}$. 
Manifold $M_3$ is Montesinos's quaternoinic space
\cite[page 120]{Montesinos87} while $M_6$ is his 
octahedral space \cite[page 117]{Montesinos87}. 
This leaves the issue of whether $M_2$ and $M_5$ are isometric, for they
have the same homology. 
Now, $M_2$ arises from a subgroup $K_2$ of the group $\Gamma$ with symbol,
$$
\begin{picture}(76,4)
\put(12,5){$4$}\put(36,5){$$}\put(60,5){$$}
\put(2,2){\circle{4}}\put(4,2){\line(1,0){20}}
\put(26,2){\circle{4}}\put(28,2){\line(1,0){20}}
\put(50,2){\circle{4}}\put(52,2){\line(1,0){20}}
\put(74,2){\circle{4}}\end{picture}.
$$
By Theorem \ref{thm2}(1), the order of $K_2$
is the index in $\Gamma$ of $\Gamma_v$, and since $|\Gamma|=2^44!$ (it
is the Weyl group $B_4$) and $|\Gamma_v|=48$, the number of symmetries
of a cube, we have $|K|=8$ (in fact it turns out that
$\pi_1(M_2)\cong\Z_8$ with generator $x_3x_2x_1x_4$,
where $x_i$ is the generator of $\Gamma$ corresponding to the
$i$-th node from the left in the symbol).
On the other hand, by the same argument, the group $K_5$
yielding $M_5$ must have order $24$ ($\Gamma$ in this case is the Weyl group $F_4$
of order $1152$). Thus, the two fundamental groups are not isomorphic,
and so the manifolds are non-isometric.

Table 2 gives the Euclidean manifolds, with $M_{10}$ the 3-torus. 
Unfortunately, we were not able to determine, by the techniques of this
paper, whether $M_8$ and $M_{12}$ were isometric or distinct.\footnote{I
Prok \cite{Prok96} has shown that $M_8$ and $M_{12}$ are related by a
Euclidean similarity and so are indeed the same manifold.}

Table 3 gives the hyperbolic results. Manifolds $M_{13}$ and $M_{14}$ come from the octahedron in
$\{4,4,3\}$, $M_{15}$, $M_{16}$ and $M_{17}$ from the cube in $\{4,3,6\}$ (see also \cite{Prok96}) and $M_i$, 
$i=18$ to $27$, from
the dodecahedron in $\{5,3,6\}$. 
Manifold $M_{14}$ is the Whitehead link complement \cite[Section
3.3]{Thurston97}.
The tetrahedron in $\{3,3,6\}$ gave no orientable manifolds, although
the non-orientable Gieseking manifold of 1911 is known to arise from it.



Manifolds $M_{13}$ and $M_{14}$ are non-isometric, despite having the same
first homology, for, using low index subgroups in Magma again, $K_{13}$ has five conjugacy classes of index 3 subgroups
while $K_{14}$ has six, so these two groups cannot be conjugate.
For the same reason, $M_{15}$ and $M_{17}$ are distinct. 
Now the group $\Gamma=\{4,3,6\}$ is arithmetic by \cite{Vinberg67}, and
thus the subgroups $K_{15}$ and $K_{17}$ are too. On the otherhand, by
the comments at the end of Section \ref{sect3}, $K_{19},K_{20}$ and $K_{26}$
are non-arithmetic, so cannot be isomorphic to $K_{15}$ and $K_{17}$.
Hence $M_{15}$ and $M_{17}$ are not isometric to any of $M_{19},M_{20}$ or
$M_{26}$. 

Finally, there are a number of pairs with the same first homology among
the $M_i$ for $i=18\text{ to }27$. Clearly $M_{22}$ and $M_{24}$ must be distinct,
for they have a different number of cusps. In fact, all ten are distinct:
the corresponding $K_i$ are non-conjugate in $\{5,3,6\}$ by
construction, and then Theorem \ref{margulis} and the comments at the end of
Section \ref{sect3} give that they are non-conjugate in $G=PO_{1,3}(\R)$. 



\small{}


\begin{thebibliography}{9}

\bibitem{Adams92}
C C Adams.
\newblock Noncompact hyperbolic 3-orbifolds of small volume.
\newblock in {\em Topology '90\/} (Columbus, Ohio 1990) 1--15, Ohio State Univ.
Math. Res. Inst. Publ. 1, de Gruyter, Berlin, 1992.

\bibitem{Best71}
L A Best.
\newblock On Torsion-free discrete subgroups of
$\psl_2({\ams C})$ with compact orbit space.
\newblock {\em Can. J. Math.}, 23 (1971), 451--460.

\bibitem{Borel81}
A Borel.
\newblock Commensurability classes and volumes of hyperbolic 3-manifolds.
\newblock {\em Ann. Sc. Norm. Sup. Pisa\/}, 8 (1981), 1-33.

\bibitem{Borel62}
A Borel and Harish-Chandra.
\newblock Arithmetic subgroups of algebraic groups.
\newblock {\em Ann. of Math.\/}, 75(3) (1962), 485-535.

\bibitem{Bourbaki68}
N Bourbaki.
\newblock {\em Groupes et alg\'{e}bres de Lie}.
\newblock Chapters 4-6, Hermann, Paris 1968; Masson, Paris 1981.

\bibitem{Cannon84}
J J Cannon.
\newblock {\em An introduction to the group theory
language Cayley}.
\newblock Academic Press, San Diego, 1984.

\bibitem{Costa92}
A F Costa.
\newblock Locally regular coloured graphs.
\newblock {\em J. Geometry}, 43 (1992), 57--74.


\bibitem{Hilbert52}
D Hilbert and S Cohn-Vossen.
\newblock {\em Geometry and the
Imagination}.
\newblock Chelsea Publishing Co., New York 1952.

\bibitem{Howlett93}
B Brink and R Howlett.
\newblock A finiteness property and an automatic structure for Coxeter groups.
\newblock {\em Math. Ann.\/} 296 (1993), 179-190.

\bibitem{Humphreys90}
J E Humphreys.
\newblock {\em Reflection groups and Coxeter groups}.
\newblock Cambridge Advanced studies in Mathematics 29, CUP 1990.

\bibitem{Lorimer92}
P J Lorimer.
\newblock Four Dodecahedral spaces.
\newblock {\em Pac. J. Math.}, 156 2 (1992), 329--335.

\bibitem{Molnar92}
E Molna'r.
\newblock Polyhedron complexes with simply transitive group actions 
and their realizations.
\newblock {\em Acta. Math. Hungar.}, 59 (1992), 175--216.

\bibitem{Montesinos87}
J M Montesinos.
\newblock {\em Classical Tessellations and
Three-Manifolds}.
\newblock Universitext, Springer 1987.


\bibitem{Prok98}
I Prok.
\newblock Classification of dodecahedral space forms.
\newblock {\em Beitr\"{a}ge Algebra Geom.}, 39(2) (1998), 497--515.

\bibitem{Prok96}
I Prok.
\newblock Fundamental tilings with marked cubes in spaces of constant curvature.
\newblock {\em Acta. Math. Hungar.}, 71 (1996), 1--14.

\bibitem{Ratcliffe94}
J G Ratcliffe.
\newblock {\em Foundations of hyperbolic manifolds}.
\newblock Graduate Texts in Mathematics 149, Springer 1994.


\bibitem{Rubinstein82}
J  Richardson and J H Rubinstein.
\newblock Hyperbolic manifolds from a regular polyhedron.
\newblock preprint.

\bibitem{Thurston97}
W Thurston.
\newblock {\em Three dimensional Topology and Geometry}.
\newblock Princeton Mathematical Series, 1997.

\bibitem{Vinberg67}
E B Vinberg.
\newblock Discrete groups in Lobachevskii spaces generated by reflections
\newblock {\em Math. Sb.\/} 72 (1967), 471-488.
\newblock $=$ {\em Math USSR-Sb.\/} 1 (1967), 429-444.

\bibitem{Weber33}
C Weber and H Seifert.
\newblock Die Beiden
Dodekaeder{\"{a}}ume.
\newblock {\em Math. Z.}, 37 (1933), 237-253.

\end{thebibliography}
\end{document}